\documentclass{amsart}

\newtheorem{theorem}{Theorem}[section]

\newtheorem{lemma}[theorem]{Lemma}
\newtheorem{proposition}[theorem]{Proposition}
\newtheorem{corollary}[theorem]{Corollary}

\theoremstyle{definition}
\newtheorem{definition}[theorem]{Definition}

\newtheorem{remark}[theorem]{Remark}

\theoremstyle{remark}

\DeclareFontFamily{U}{wncy}{}
\DeclareFontShape{U}{wncy}{m}{n}{<->wncyr10}{}
\DeclareSymbolFont{mcy}{U}{wncy}{m}{n}
\DeclareMathSymbol{\Sh}{\mathord}{mcy}{"58}

\newcommand\mylabel[1]{\label{#1}}

\newcommand{\ZZ}{\mathbb{Z}}

\newcommand{\GG}{\mathbb{G}}

\newcommand  {\shA}     {\mathcal{A}}

\newcommand  {\shC}     {\mathcal{C}}

\newcommand  {\shEnd}   {\mathcal{E}\!\text{\textit{nd}}}

\newcommand  {\shD}     {\mathcal{D}}
\newcommand  {\shE}     {\mathcal{E}}
\newcommand  {\shF}     {\mathcal{F}}
\newcommand  {\shG}     {\mathcal{G}}
\newcommand  {\shH}     {\mathcal{H}}
\newcommand  {\shHom}   {\mathcal{H}\!\text{\textit{om}}}
\newcommand  {\shI}     {\mathcal{I}}

\newcommand  {\shK}     {\mathcal{K}}
\newcommand  {\shM}     {\mathcal{M}}

\newcommand  {\shL}     {\mathcal{L}}

\def\stA{\mathscr{A}}

\def\stG{\mathscr{G}}

\def\stX{\mathscr{X}}

\newcommand  {\Ab}      {{\text{\rm Ab}}}

\newcommand  {\Aff}     {\text{Aff}}

\newcommand  {\Aut}     {\operatorname{Aut}}

\newcommand  {\biget}   {{\text{\rm big-et}}}
\newcommand  {\Biget}   {\operatorname{Big-et}}
\newcommand  {\Br}      {\operatorname{Br}}
\newcommand  {\BBr}     {\widetilde{\operatorname{Br}}}

\newcommand  {\GL}      {\operatorname{GL}}

\newcommand  {\id}      {\operatorname{id}}

\newcommand  {\Liset}   {\operatorname{Lis-et}}
\newcommand  {\liset}   {{\text{\rm lis-et}}}

\newcommand  {\dirlim}  {\varinjlim}

\newcommand  {\lra}     {\longrightarrow}

\newcommand  {\Mor}     {\operatorname{Mor}}

\renewcommand{\O}       {\mathcal{O}}

\newcommand  {\pr}      {\operatorname{pr}}

\newcommand  {\ra}      {\rightarrow}

\newcommand  {\Spec}    {\operatorname{Spec}}

\newcommand  {\Split}   {\operatorname{Split}}

\def\mydate{\number\day\space\ifcase\month \or January\or February\or March\or 
April\or May\or June\or July\or
August\or September\or October\or November\or December\fi \space\number\year}

\DeclareFontFamily{U}{wncy}{}
\DeclareFontShape{U}{wncy}{m}{n}{<->wncyr10}{}
\DeclareSymbolFont{mcy}{U}{wncy}{m}{n}
\DeclareMathSymbol{\Sh}{\mathord}{mcy}{"58}

\usepackage{amscd,amssymb,amsmath}
\usepackage{mathrsfs}
\usepackage[all]{xy}

\begin{document}

\title[Twisted sheaves]
      {The bigger Brauer group and twisted sheaves}

\author[Jochen Heinloth]{Jochen Heinloth}
\address{
Korteweg-de Vries Institute for Mathematics,
 University of Amsterdam, 
1018 TV Amsterdam, 
The Netherlands}
\curraddr{}
\email{heinloth@science.uva.nl}

\author[Stefan Schr\"oer]{Stefan Schr\"oer}
\address{Mathematisches Institut, 
Heinrich-Heine-Universit\"at,
40225 D\"usseldorf, 
Germany}
\curraddr{}
\email{schroeer@math.uni-duesseldorf.de}

\subjclass[2000]{14A20, 14F22}

\dedicatory{25 March 2007}

\begin{abstract}
Given an algebraic stack with quasiaffine diagonal,
we show that each $\GG_m$-gerbe comes from 
a central separable algebra. 
In other words,   Taylor's bigger Brauer group
equals the   \'etale cohomology in degree two with coefficients
in $\GG_m$. This gives new results also for schemes.
We use the method of twisted sheaves
explored by de Jong and Lieblich.
\end{abstract}

\maketitle
\tableofcontents

\section*{Introduction}
\mylabel{introduction}

Let $X$ be a scheme.
About forty years ago, Grothendieck \cite{GB} posed the
problem wether  the inclusion $\Br(X)\subset H^2(X,\GG_m)$
of the Brauer group 
of Azumaya algebras 
coincides with the torsion part
of the \'etale cohomology group. 
It is known that this fails for certain nonseparated schemes
\cite{Edidin; Hassett; Kresch; Vistoli 1999}.
On the other hand, there are strong positive results.
Gabber \cite{GB} proved equality  for affine schemes,
and also had an unpublished proof   for schemes carrying  ample line bundles.
Recently, de Jong \cite{de Jong 2006} gave a new proof for the latter statement,
 based on the notion of 
\emph{twisted sheaves},  that is, sheaves
on gerbes. This method turn out to be rich and powerful,
and was further explored by Lieblich in \cite{Lieblich 2005} and
\cite{Lieblich 2007}. 

In this paper we shall prove  that there is, for arbitrary noetherian schemes,
an equality $\BBr(X)=H^2(X,\GG_m)$, where
$\BBr(X)$ is the \emph{bigger Brauer group}.
This group is defined in terms of so-called
\emph{central separable algebras}, and was introduced
by Taylor \cite{Taylor 1982} (Caenepeel and Grandjean  \cite{Caenepeel; Grandjean 1998} later fixed some technical problem
in the original definition).
Such algebras are defined and behave very similar to Azumaya algebras,
but do not necessarily contain a unit.
Raeburn and Taylor  \cite{Raeburn; Taylor 1985}  constructed an inclusion $\BBr(X)\subset H^2(X,\GG_m)$ using
methods from nonabelian cohomology, and   showed that
this inclusion actually an equality   provided $X$ 
carries   ample line bundles.
To remove this assumption, we shall use de Jong's insight \cite{de Jong 2006}
and work with a gerbe $\stG$ defining the cohomology class $\alpha\in H^2(X,\GG_m)$.
The basic observation   is that $\stG$ may be  viewed as
an algebraic stack (= Artin stack),
and that the existence of the desired central separable algebra on $X$
is equivalent to the existence of certain coherent sheaves on $\stG$.
A key ingredient in our arguments is the result of Laumon and Moret-Bailly
that quasicoherent sheaves on noetherian algebraic stacks are
direct limits of coherent sheaves \cite{Laumon; Moret-Bailly 2000}.

This stack-theoretic approach suggests a generalization of the problem at hand: 
Why not replace the scheme $X$ by an algebraic stack $\stX$?
Our investigation actually takes place in the setting.
Here, however, one has to make an additional assumption. Our main result is that
$\BBr(\stX)= H^2(\stX,\GG_m)$ holds for noetherian algebraic stacks whose  diagonal morphism $\stX\ra\stX\times\stX$
is quasiaffine.  Deligne--Mumford stacks, and in particular  
algebraic spaces and schemes, automatically satisfy this assumption. 
In contrast,  there are algebraic stacks
 with $\BBr(\stX)\subsetneq H^2(\stX,\GG_m)$.
We give an example based on   observations of Totaro \cite{Totaro 2004}.

Working with sheaves and cohomology on algebraic stacks $\stX$
requires some care.
A convenient setting is the so-called \emph{lisse-\'etale site} $\Liset(\stX)$.
For our purposes, it is useful to work a larger site as well,
which we call the \emph{big-\'etale site} $\Biget(\stX)$.
The relation between the associated topoi $\stX_\liset$ and $\stX_\biget$
is not so straightforward as one might expect at first glance.
The problem is, roughly speaking, that they are not related by 
a continuous map.
Such phenomena gained notoriety in the theory of algebraic stacks, 
and were explored 
by Behrend \cite{Behrend 2003} and Olsson \cite{Olsson 2007}.
However, in the Appendix we show that an abelian big-\'etale sheaf and its
restriction to the lisse-\'etale site have the same cohomology,
by reexamining Grothendieck's original construction of injective objects via
transfinite induction. This result appears to be of independent interest.

\section{Gerbes on algebraic  stacks}
\mylabel{algebraic stacks}

In this section we recall some   basic facts on gerbes over algebraic stacks.
Throughout, we closely follow the book of Laumon and Moret-Bailly \cite{Laumon; Moret-Bailly 2000}
in terminology and notation.

Fix a base scheme $S$, and let $(\Aff/S)$ be the category of affine schemes endowed
with a morphism to $S$. Let $\stX$ be an algebraic $S$-stack.
A \emph{lisse-\'etale sheaf} on $\stX$ is, by definition, a sheaf on the
\emph{lisse-\'etale site} $\Liset(\stX)$.
The objects of the latter are pairs $(U,u)$, where $U$ is an algebraic space,
and $u:U\ra\stX$ is a smooth morphisms.
The morphisms from $(U_1,u_1)$ to another object $(U_2,u_2)$ 
are pairs $(f,\alpha)$, where $f:U_1\ra U_2$ is a morphism of algebraic spaces,
and $\alpha$ is a natural transformation between the functors 
$u_1,u_2\circ f:U_1\ra\stX$,
such that we have a 2-commutative diagram
$$
\xymatrix{
 U_1  \ar@/^1.3pc/[rr]^{u_1}="1"  \ar@/_1.3pc/[rr]_{u_2\circ f}="2"
 && \stX.  
\ar@{}"1";"2"|(.3){\,}="7"
\ar@{}"1";"2"|(.7){\,}="8"
\ar@{=>}"7" ;"8"^{\alpha} } 
$$
The Grothendieck topology is generated by   those   $(f,\alpha)$ with   $f:U_1\ra U_2$ \'etale and surjective.
We denote by $\stX_{\liset}$ the associated \emph{lisse-\'etale  topos},
that is, the category of lisse-\'etale sheaves on $\stX$.

For the applications we have in mind, it is  natural to work with a larger site as well.
It resembles the big site of a topological space, so we call it
the \emph{big-\'etale site} $\Biget(\stX)$.
Here the objects are pairs $(U,u)$, where again $U$ is an algebraic space,
but now the morphism $u:U\ra\stX$ is arbitrary.
Morphisms and   Grothendieck topology are defined as for the lisse-\'etale site.
The associated big-\'etale topos is denoted by $\stX_{\biget}$.
Following \cite{Laumon; Moret-Bailly 2000}, we write
$\underline{\Biget}(\stX)\subset\Biget(\stX)$ for
the subsites of objects $(U,u)$ with $U$ affine.
According to the Comparison Lemma \cite{SGA 4b}, Expos\'e III, Theorem 4.1,
 this inclusion induces an
equivalence between the corresponding topoi of sheaves.

The relation between the big-\'etale and the lisse-\'etale site is not as
straightforward as one might expect.
This is because the inclusion   $\Liset(\stX)\subset\Biget(\stX)$
does not induce a map of topoi, as discussed in the Appendix.
However, given a big-\'etale sheaf $\shF$, there actually is a canonical map
$$
H^i(\stX_\biget,\shF)\lra H^i(\stX_\liset,\shF|_{\Liset(\stX)}),
$$
and we shall prove in the appendix  that this map is   bijective
(Theorem \ref{cohomology bijective}).
We therefore write $H^i(\stX,\shF)$ for the cohomology of   big-\'etale sheaves, provided there is no risk of confusion.

Next we recall some facts on  gerbes on $\stX$.
Let $G\ra S$ be an abelian group algebraic space over $S$.
(For our applications we merely use the case   $G=\GG_{m,S}$.)
This yields an abelian big-\'etale sheaf denoted $G_\stX$,
whose groups of local sections $\Gamma((U,u),G_\stX)$ is the set of $S$-morphisms
$g:U\ra G$.
In turn, we  have cohomology groups $H^i(\stX,G_\stX)$.
According to the previous paragraph, it does matter wether we compute cohomology on the lisse-\'etale or big-\'etale site.

As explained in the Giraud's treatise \cite{Giraud 1971}, Chapter IV, \S 3.4,
cohomology classes from $H^2(\stX, G_\stX)$ 
correspond to equivalence classes
of $G_\stX$-gerbes $\stG\ra\underline{\Biget}(\stX)$; equivalently, a gerbe on the lisse-\'etale site.
By composing with $(U,u)\mapsto U$, we obtain a functor $\stG\ra(\Aff/S)$.
As  Lieblich explains in \cite{Lieblich 2007}, Proposition 2.4.3, 
this makes $\stG$ into an $S$-stack, endowed with a 1-morphism of $S$-stacks
$F:\stG\ra\stX$. Under fairly general assumptions,
this $S$-stack   is algebraic; the following criterion 
generalizes a result of de Jong \cite{de Jong 2006} and  Lieblich (\cite{Lieblich 2007}, Corollary 2.4.4):

\begin{proposition}
\mylabel{gerbe algebraic}
Notation as above.
Suppose that the structure morphism $G\ra S$ is smooth and of finite presentation.
Then the $S$-stack $\stG$ is algebraic. 
\end{proposition}

\proof
Choose a smooth surjection $P:X\ra\stX$ for some scheme $X$,
with $\stG_{X,P}$   nonempty. Then
the projection $\stG\times_\stX X\ra X$ has a section, and by
 \cite{Laumon; Moret-Bailly 2000}, Lemma 3.21,
there is a 1-isomorphism $\stG\times_\stX X\ra B(G_X)$
into the $S$-stack of $G_X$-torsors and this stack is algebraic.

Choose a smooth surjection $Y\ra\stG\times_\stX X$ from
some   scheme $Y$.
Composing with the second projection, we obtain a smooth, surjective,
representable morphism $H:Y\ra\stG$.
In light of loc.\ cit., Proposition 4.3.2,
it remains to check that the canonical morphism   
$Y\times_\stG Y\ra Y\times Y$ is quasicompact and separated.
Note that both $S$-stacks in question are associated to schemes.
Our  map factors over $Y\times_\stX Y$, and
the morphism of schemes $Y\times_\stX Y\ra Y\times Y$ is   quasicompact and separated, because the $S$-stack $\stX$ is algebraic.
Whence it suffices to check that $Y\times_\stG Y\ra Y\times_\stX Y$
is quasicompact and separated.

To verify this, consider the following  $G$-action on the objects of the
$S$-stack $Y\times_\stG Y$:
Given some $U\in(\Aff/S)$, 
the objects in $Y\times_\stG Y$ over $U$ are, by definition,
triples $(u_1,u_2,\varphi)$, where $u_i:U\ra Y$ are $S$-morphisms,
and $\varphi:H(u_1)\ra H(u_2)$ is an isomorphism in $\stG_U$.
Then the $S$-morphisms $g:U\ra G$ act on such tripels via
$
g\cdot(u_1,u_2,\varphi) = (u_1,u_2,g\varphi).
$
Using   that $\stG\ra\underline{\Biget}(\stX)$ is a $G_\stX$-gerbe,
we infer that our morphism   $Y\times_\stG Y\ra Y\times_\stX Y$,
viewed as a morphism of schemes, is a $G$-principle bundle with respect to the \'etale topology.
Our assumptions on the structure morphism $G\ra S$ ensure
that it is quasicompact and separated.
By descent, the  same   holds for $Y\times_\stG Y\ra Y\times_\stX Y$, see
\cite{SGA 1}, Expos\'e V, Corollary 4.6 and 4.8.
\qed

\begin{remark}
 Using Artin's theorem  \cite{Laumon; Moret-Bailly 2000} Proposition 10.31.1, the above proof generalizes to the case that $\stG\to S$ flat group schemes of finite presentation if one considers gerbes in the fppf-topology.
\end{remark}

\medskip

Since we assumed the structure morphism $G\ra S$ to be smooth, 
it is easy to see that the resulting morphism $F:\stG\ra\stX$ of
algebraic $S$-stacks is smooth as well, compare \cite{Laumon; Moret-Bailly 2000},
Remark 10.13.2. Given a  quasicoherent sheaf $\shH$ on $\stX$,
we obtain functorially a quasicoherent  sheaf $F^*(\shH)$ on $\stG$,
defined by 
$$
F^*(\shH)_{U,u}=\shH_{U,Fu},\quad (U,u)\in\Liset(\stG).
$$
We now describe those quasicoherent sheaves on $\stG$ that are of isomorphic
to pullbacks $F^*(\shH)$.
Let $\shF$ be a quasicoherent sheaf on $\stG$, and $(U,u)\in\Liset(\stG)$.
Any local section  $g\in\Gamma((U,Fu),G_\stX)$ induces an automorphism
$(\id_U,g):(U,u)\ra (U,u)$ in the lisse-\'etale site.
In turn, it acts bijectively   on   local sections
\begin{equation}
\label{induced bijection}
(\id_U,g)^*:\Gamma((U,u),\shF)\lra\Gamma((U,u),\shF).
\end{equation}
Sheaves for which all these bijections are actually identities shall play
an important role throughout. Let us introduce the following terminology,
which comes from the special case $G=\GG_{m,S}$:

\begin{definition}
\mylabel{weight zero}
A quasicoherent sheaf $\shF$ on $\stG$ is called \emph{of weight zero}
if the bijections in (\ref{induced bijection})
are identities for all $(U,u)$ and $g$ as above.
\end{definition}

The following characterization of sheaves of weight zero is well-known:

\begin{lemma}
\mylabel{descent}
The functor $\shH\mapsto F^*(\shH)$ is an equivalence between the
category of quasicoherent sheaves on $\stX$ and the category of
quasicoherent sheaves on $\stG$ of weight zero.
\end{lemma}

\proof
Choose a smooth surjection $u:U\ra\stG$ from some scheme $U$.
According to \cite{Laumon; Moret-Bailly 2000}, Proposition 13.2.4,
the category of quasicoherent sheaves  on $\stG$ is equivalent
to the category of quasicoherent sheaves on $U$ endowed with
a descent datum with respect to $u$.
Let $\shF$ be a quasicoherent sheaf on $\stG$ of weight zero,
with induced descent datum $\varphi:\pr_1^*(\shF_{U,u})\ra\pr_2^*(\shF_{U,u})$
on $U\times_\stG U$. As discussed in the proof of Proposition \ref{gerbe algebraic},
the morphism $U\times_\stG U\ra U\times_\stX U$
 is a $G_{U\times_\stX U}$-torsor.  Since $\shF$ is of weight zero,
$\varphi$ is invariant under $G_{U\times_\stX U}$, whence
descends to $ U\times_\stX U$. In this way we obtain for the quasicoherent
sheaf $\shF_{U,u}$ on $U$
a descent datum with respect to the smooth surjection $Fu:U\ra\stX$,
which in turn defines a quasicoherent sheaf $\shH$ on $\stG$.
It is easy to see that there is a natural isomorphism $\shF\simeq F^*(\shH)$,
and that the functor $\shF\mapsto \shH$ is quasiinverse to $\shH\mapsto F^*(\shH)$.
\qed

\section{Taylor's bigger Brauer group}
\mylabel{brauer group}

In this section we recall and discuss Taylor's  
bigger Brauer group  \cite{Taylor 1982} in the general
context of algebraic stacks.
Taylor's   idea is to attach to certain kinds of (not necessarily unital) associative algebras on $\stX$ a
$\GG_m$-gerbe, which in turn yields a cohomology class
in $H^2(\stX,\GG_m)$. The collection of all such 
cohomology classes constitutes a subgroup, which is called
the \emph{bigger Brauer group} $\widetilde{\Br}(\stX)\subset H^2(\stX,\GG_m)$.

Let us now go into details.
Suppose we have two quasicoherent sheaves $\shM$ and $\shH$ on $\stX$,
together with a pairing $\Phi:\shH\otimes\shM\ra\O_\stX$.
This defines a quasicoherent  associative $\O_\stX$-algebra 
$\shM\otimes^\Phi\shH$ as follows: The underlying  quasicoherent sheaf is $\shM\otimes\shH$;
the multiplication law is defined on local sections by
$$
(m\otimes h)\cdot (m'\otimes h')= m\otimes\Phi(h,m') h'.
$$
An important special case is that $\shM$ is locally free of finite rank,
$\shH=\shM^\vee$ is the dual sheaf, and $\Phi(h,m)=h(m)$ is the evaluation pairing.
Then   $\shM\otimes^\Phi\shH$ is canonically isomorphic to the endomorphism
algebra $\shEnd(\shM)$, which contains a unit. Note, however, that in general  
$\shM\otimes^\Phi\shH$ does  not   contain a unit.

In the following we are interested in  algebras    that
are locally of  the  form $\shM\otimes^\Phi\shH$, where one additionally demands that the pairing
$\Phi:\shH\otimes\shM\ra\O_\stX$ is surjective. 
Given an $\O_\stX$-algebra $\shA$, we use the following ad hoc terminology: A \emph{local splitting}
for $\shA$ is a sextuple $(U,u,\shM,\shH,\Phi,\psi)$,
where $U$ is an algebraic space, $u:U\ra\stX$ is a morphism of $S$-stacks,
$\shM$ and $\shH$ are quasicoherent $\O_U$-modules,
$\Phi:\shH\otimes\shM\ra\O_U$ is a surjective linear map,
and $\psi:\shM\otimes^\Phi\shH\ra\shA_{U,u}$ is an bijection of algebras.

The local splittings form a category: A \emph{morphism between  two local splittings}
$(U,u,\shM,\shH,\Phi,\psi)$ and $(U',u',\shM',\shH',\Phi',\psi')$
is a quadruple $(f,\alpha,s,t)$, where $(f,\alpha)$
is a morphism from $u:U\ra\stX$ to $u':U'\ra\stX$, and
$s:f_*(\shM)\ra\shM'$ and $t:f_*(\shH)\ra\shH'$ are linear maps of  
sheaves on $U'$; we demand that the adjoint maps $\shM\ra f^*(\shM')$ and $\shH\ra f^*(\shH')$ are bijective and  that the diagram 
$$
\begin{CD}
f_*(\shM\otimes^\Phi\shH) @>s\otimes t>> \shM'\otimes^{\Phi'}\shH'\\
@V\psi VV @VV\psi'V\\
f_*(\shA_{U,u}) @>>\text{can} > \shA_{U',u'}
\end{CD}
$$
is commutative. Composition is defined in the obvious way.

Let $\Split(\shA)$ denote the category of all local splittings of $\shA$ with $U$ affine.
Then we have a forgetful functor
$$
\Split(\shA)\lra(\Aff/S),\quad (U,u,\shM,\shH,\Phi,\psi)\longmapsto U,
$$
which gives  $\Split(\shA)$ the structure of an  $S$-stack.
It comes along with a 1-morphism of $S$-stacks $\Split(\shA)\ra\stG$,
sending a local splitting $(U,u,\shM,\shH,\Phi,\psi)$ to the object in $\stG_{U,u}$ induced by the morphism
$u:U\ra\stG$.
Moreover, $(U,u,\shM,\shH,\Phi,\psi)\mapsto(U,u)$ makes $\Split(\shA)$ into a stack
over the site $\underline{\Biget}(\stX)$.

Given a local  section $s\in\Gamma((U,u),\GG_{m,\stX})=\Gamma(U,\O_U^\times)$
and a local splitting $(U,u,\shM,\shH,\Phi,\psi)\in\stG_{U,u}$,
we obtain an automorphism $(\id_U,\id_u,s,s^{-1})$ on this object.
According to the result of Raeburn and Taylor (\cite{Raeburn; Taylor 1985}, Lemma 2.4) 
the resulting map of sheaves
$$
\O_\stX^\times|_{(\Aff/U)}\lra\underline{\Aut}_{\stG}(U,u,\shM,\shH,\Phi,\psi)
$$
is bijective; moreover,   all objects from $\stG_{U,u}$ are locally isomorphic.
So if we demand that the algebra $\stA$ on $\stX$ admits a splitting
over some $u:U\ra\stX$ that is smooth and surjective,
the   stack $\stG\ra \underline{\Biget}(\stX)$ is   a $\GG_{m,\stX}$-gerbe, whence
yields a cohomology class $[\shA]\in H^2(\stX,\GG_m)$:

\begin{definition}
\mylabel{}
The algebra $\shA$ on $\stX$
is called a \emph{central separable algebra} if it admits a local splitting
$(U,u,\shM,\shH,\Phi,\psi)$ with  $u:U\ra\stX$   smooth surjective.
\end{definition}

Note that this   
differs slightly from Taylor's approach in  \cite{Taylor 1982}, Definition 2.1.
By taking the existence of splittings as defining property, and not
as a consequence,
we  avoid the technical problems discussed 
in \cite{Caenepeel; Grandjean 1998}.

We define the  \emph{bigger Brauer group} $\BBr(\stX)\subset H^2(\stX,\GG_m)$
as the subgroup generated by cohomology classes coming from central separable algebras as described above.
Our task is to find conditions implying that the inclusion
$\BBr(\stX)\subset H^2(\stX,\GG_m)$ is actually an equality. The following 
properties of quasicoherent sheaves will be useful:

\begin{proposition}
\mylabel{equivalent conditions}
Let $\shF$ be a quasicoherent sheaf on an algebraic $S$-stack $\stG$.
The following two conditions are equivalent:
\renewcommand{\labelenumi}{(\roman{enumi})}
\begin{enumerate}
\item There is a smooth surjection $u:U\ra\stG$ from an algebraic space $U$
and a surjective linear map $\shF_{U,u}\ra\O_U$.
\item
There is a smooth surjection $v:V\ra\stX$ from an algebraic space $V$
and a   decomposition $\shF_{V,v}\simeq\shK\oplus\O_V$ for some
quasicoherent sheaf $\shK$ on $V$.
\end{enumerate}
\end{proposition}

\proof
The implication (ii)$\Rightarrow$(i) is trivial. To see 
(i)$\Rightarrow$(ii), suppose we have a surjection $\shF_{U,u}\ra\O_U$.
Choose an \'etale surjection $V\ra U$, where $V=\bigcup V_\alpha$ 
is a disjoint union of affine schemes.
Let $v:V\ra\stG$ be the induced morphism, and $\shK$ be the kernel of the
induced surjection $\shF_{V,v}\ra\O_V$. This surjection must have a section,
because quasicoherent sheaves on affine schemes have no higher cohomology.
\qed

\medskip
Let us introduce a name for such sheaves:

\begin{definition}
\mylabel{invertible summand}
Let $\shF$ be a quasicoherent sheaf on an algebraic $S$-stack $\stG$.
We say that $\shF$ \emph{locally contains invertible summands} if
the two equivalent conditions of Proposition \ref{equivalent conditions} hold.
\end{definition}

This notion was used in \cite{Schroeer 2007}
to solve some problems on singularities in positive characteristics.
For coherent sheaves  on noetherian stacks, 
we have the following characterization involving the dual 
sheaf $\shF^\vee=\shHom(\shF,\O_\stG)$:

\begin{proposition}
\mylabel{evaluation pairing}
Let $\shF$ be a coherent sheaf on a noetherian algebraic $S$-stack $\stG$.
Then the following are equivalent:
\renewcommand{\labelenumi}{(\roman{enumi})}
\begin{enumerate}
\item The   sheaf $\shF$ locally contains invertible direct summands.
\item The evaluation pairing $\shF\otimes\shF^\vee\ra\O_\stG$ is surjective.
\item There is a smooth surjective morphism $u:U\ra\stG$ from some affine scheme 
$U=\Spec(R)$, an $R$-module $N$, and  a   surjective linear mapping $\Gamma((U,u),\shF)\otimes_R N\ra R$.
\end{enumerate}
\end{proposition}

\proof
The implication (i)$\Rightarrow$(ii) is trivial:
Choose a smooth surjection $u:U\ra\stG$ from some affine scheme $U$
so that  $\shF_{U,u}\simeq\shK\oplus\O_U$.
Then the evaluation paring $\shF_{U,u}\oplus\shF_{U,u}^\vee\ra\O_U$
is obviously surjective, and so is $\shF\otimes\shF^\vee\ra\O_\stG$.
The implication (ii)$\Rightarrow$(iii) is also trivial:
Choose any smooth surjection $u:U\ra\stG$ from some
affine scheme $U$ and set $N=\Gamma((U,u),\shF^\vee)$.

It remains to check (iii)$\Rightarrow$(i).
Choose a smooth surjection  $u:U\ra\stG$   from some   
affine scheme $U=\Spec(R)$ admitting a surjection
$\phi:\Gamma((U,u),\shF)\otimes_R N\ra R$.
Then there are finitely many $f_1,\ldots, f_r\in \Gamma((U,u),\shF)$
and $n_1,\ldots,n_r\in N$ with $\phi(\sum f_i\otimes n_i)=1$.
Setting $s_i=\phi(f_i\otimes n_i)$, we obtain
an affine open covering $U=V(s_1)\cup\ldots\cup V(s_r)$.
Replacing $U$ by the disjoint union of the $V(s_i)$, we easily
reduce to the case $r=1$. This means that there is
an $f\in\Gamma((U,u),\shF)$ and $n\in N$ with $\varphi(f\otimes n)=1$.
In other words, the map $f\mapsto \phi(f\otimes n)$ is surjective,
which gives the desired surjection $\shF_{U,u}\ra\O_U$.
\qed

\medskip
We finally examine the connection  to  central separable algebras.
Suppose $\stX$ is an algebraic $S$-stack, and $\stG\ra\underline{\Biget}(\stX)$
is a $\GG_{m,\stX}$-gerbe.
Let $F:\stG\ra\stX$ be the resulting morphism of algebraic $S$-stacks, as discussed
in Section \ref{algebraic stacks}.
Given a   lisse-\'etale sheaf $\shF$ on $\stG$ and a smooth morphism
$u:U\ra\stG$ from some algebraic space $U$, we denote
by $\shF_{U,u}$ the induced sheaf on $U$.
For quasicoherent sheaves, the actions of $\GG_{m,U}$ on $\shF_{U,u}$
corresponds to a weight decomposition $\shF=\bigoplus\shF_n$,
as explained in \cite{SGA 3a}, Expos\'e I, Proposition 4.7.2.
Here the direct sum runs through all $n\in\ZZ$, which is the character group
of $\GG_m$.
A quasicoherent sheaf with $\shF=\shF_n$ is called \emph{of weight $w=n$}.

\begin{theorem}
\mylabel{gerbe equivalent}
Let $\shG$ be a $\GG_m$-gerbe on a noetherian algebraic $S$-stack $\stX$.
Then the following are equivalent:
\renewcommand{\labelenumi}{(\roman{enumi})}
\begin{enumerate}
\item
There is a central separable algebra $\shA'$ on $\stX$ whose $\GG_m$-gerbe of splittings $\Split(\shA')$ is equivalent to $\stG$.
\item
There is a \emph{coherent} central separable algebra $\shA$ on $\stX$ whose
$\GG_m$-gerbe of splittings $\Split(\shA)$ is equivalent to $\stG$.
\item
There is a coherent sheaf $\shF$ on $\stG$ of weight $w=1$ that
locally contains invertible summands.
\end{enumerate}
\end{theorem}

\proof
The implication (ii)$\Rightarrow$(i) is trivial.
To prove (i)$\Rightarrow$(iii), assume that $\stG=\Split(\shA')$ 
for some central separable algebra $\shA'$ on $\stX$.
Let $\tilde{u}:U\ra\stG$ be a smooth morphism from an affine scheme $U$,
and   $(U,u,\shM,\shH,\Phi,\psi)\in\stG_{U,\tilde{u}}$ be the resulting object,
as described in Section \ref{}.
We now use the sheaves $\shM$ on $U$ to define a sheaf $\underline{\shM}$
on $\stG$ by the tautological formula
$$
\Gamma((U,\tilde{u}),\underline{\shM})= \Gamma(U,\shM).
$$
This obviously defines a presheaf on $\stG$.
It is easy to check that
it satisfies the sheaf axiom, and that 
$\underline{\shM}_{U,\tilde{u}}\simeq\shM$,
such that $\underline{\shM}$ is quasicoherent.
This quasicoherent sheaf  is of weight $w=1$:
The sections $s\in\Gamma((U,\tilde{u}),\GG_{m,\stX})=\Gamma(U,\O_U^\times)$
act via the automorphism $(\id_U,\id_u,s,s^{-1})$ on the object $(U,u,\shM,\shH,\Phi,\psi)\in\stG_{U,u}$,
whence by   multiplication-by-$s$ on $\Gamma((U,\tilde{u}),\underline{\shM})$. 

To proceed, consider the ordered set  $\shF_\alpha\subset\underline{\shM}$,
$\alpha\in I$ of coherent subsheaves.
The induced map $\dirlim(\shF_\alpha)\ra\underline{\shM}$ is bijective,
by \cite{Laumon; Moret-Bailly 2000}, Proposition 15.4.
It remains to verify that some $\shF_\alpha$ locally contains invertible
summands.
By construction, we have $\underline{\shM}_{U,\tilde{u}}\simeq\shM$, and a surjective pairing
$\Phi:\shM\otimes\shH\ra\O_U$.
Setting $M_\alpha=\Gamma((U,\tilde{u}),\shF_\alpha)$ and $N=\Gamma(U,\shH)$,
we obtain a surjective pairing $\dirlim(M_\alpha)\otimes N\ra R$.
Using that direct limits commute with tensor products,
we infer that the map $M_\beta\otimes N\ra R$ must already by surjective for some $\beta\in I$.
According to Proposition \ref{evaluation pairing}, the sheaf $\shF=\shF_\beta$ locally contains invertible summands.

It remains to prove the implication (iii)$\Rightarrow$(ii). Let $\shF$ be a coherent sheaf
on $\stG$, of weight $w=1$ and locally containing invertible summands.
Then the evaluation paring $\Phi:\shF^\vee\otimes\shF\ra\O_\stG$ is surjective,
such that $\shF\otimes^\Phi\shF^\vee$ is a central separable algebra on $\stG$,
and the underlying coherent sheaf has weight zero.
It follows from Lemma \ref{descent} that $\shF\otimes^\Phi\shF^\vee$
is isomorphic to the preimage of a nonunital associative algebra $\shA$
on $\stX$.
Moreover, given a smooth morphism $\tilde{u}:U\ra\stG$, we easily infer
that we have a canonical isomorphism
$\psi:\shA_{U,F\tilde{u}}\ra\shF_{U,\tilde{u}}\otimes^\Phi\shF^\vee_{U,\tilde{u}}$,
whence the algebra $\shA$ is   central separable.
 
To finish the proof, we have to construct a functor of $\GG_{m,\stX}$-gerbes $\stG\ra\Split(\shA)$.
Let   $X\in\stG_{U,u}$ be an object. Choose a morphism $\tilde{u}:U\ra\stG$
inducing this object, set $\shM=\shF_{U,\tilde{u}}$ 
and $\shH=\shF^\vee_{U,\tilde{u}}$, and let $\Phi:\shM\otimes\shH\ra\O_U$
be the evaluation paring.
Together with the canonical isomorphism
 $\psi:\shA_{U,u}\ra\shF_{U,\tilde{u}}\otimes^\Phi\shF^\vee_{U,\tilde{u}}$ 
described above,
we obtain the desired functor as
$$
\stG\lra\Split(\shA),\quad X\longmapsto (U,u,\shM,\shH,\Phi,\psi),
$$
which is obviously compatible with the $\GG_{m,\stX}$-action.
\qed

\section{Existence of central separable algebras}
\mylabel{existence algebras}
 
We now come to our main result:

\begin{theorem}
\mylabel{main result}
Let $\stX$ be a noetherian algebraic $S$-stack whose diagonal morphism
$\Delta:\stX\ra\stX\times\stX$ is quasiaffine.
Then $\BBr(\stX)=H^2(\stX,\O_{\stX}^\times)$.\
\end{theorem}

Before we prove this, let us discuss two special cases.
For schemes, the diagonal morphism is an embedding, whence automatically quasiaffine.
Thus the preceding Theorems applies to schemes, which removes
superfluous assumptions in   
results of Raeburn and Taylor \cite{Raeburn; Taylor 1985} 
and the second author \cite{Schroeer 2003}.
According to \cite{Laumon; Moret-Bailly 2000}, Lemma 4.2,
the diagonal morphism is   quasiaffine even for  Deligne--Mumford stacks.
Thus:

\begin{corollary}
\mylabel{deligne-mumford stacks}
Let $\stX$ be   a noetherian scheme or a  noetherian Deligne--Mumford
$S$-stack. Then we have equality  $\widetilde{\Br}(\stX)=H^2(\stX,\GG_m)$.
\end{corollary}

\medskip
\emph{Proof of Theorem \ref{main result}:}
 Fix a cohomology class
$\alpha\in H^2(\stX,\GG_m)$   and choose a 
$\GG_m$-gerbe $\stG\ra\underline{\Biget}(\stX)$ representing $\alpha$.
Then there is an affine scheme $U$ and a smooth surjective morphism
$u:U\ra\stX$ so that $\stG_{U,u}$ is nonempty.
Note that $u:U\ra\stX$ is quasiaffine. To see this,
let $v:V\ra\stX$ be a morphism from an affine scheme $V$.
Then we have a commutative diagram with cartesian square
$$
\begin{CD}
U\times_\stX V @>>> U\times V @>\pr_2>> V\\
@VVV @VVV\\
\stX @>>\Delta> \stX\times\stX.
\end{CD}
$$
The projection $\pr_2$ is affine, because $U$ is affine.
The morphism $U\times_\stX V\ra U\times V$ is quasiaffine because $\Delta$ is quasiaffine.
Whence the composition $U\times_\stX V\ra V$ is quasiaffine, which means that  $u:U\ra\stX$ is quasiaffine. 

By assumption, the induced gerbe $\stG\times_\stX U\ra U$ is trivial.
Hence there is coherent sheaf 
$\shE$ on $\stG\times_\stX U$ of weight $w=1$  locally containing
invertible summands. Choose a smooth surjection $v:V\ra\stG\times_\stX U$  from some
affine scheme $V$ so that there is a surjection $\shE_{V,v}\ra\O_V$.

Now consider the other projection $F:\stG\times_\stX U\ra\stG$.
This morphism is quasicompact and quasiseparated, so $F_*(\shE)$ is quasicoherent.
The canonical map $F^*F_*(\shE)\ra\shE$ is surjective 
by \cite{EGA II}, Proposition 5.1.6, because $F$ is quasiaffine. Hence
the composition $F^*F_*(\shE)_{V,v}\ra\O_V$ is surjective as well.
Setting $v'=F\circ v:V\ra\stG$, we obtain a surjection $F_*(\shE)_{V,v'}\ra\O_V$.
Applying \cite{Laumon; Moret-Bailly 2000}  Proposition 15.4,
we write $F_*(\shE)=\dirlim\shF_i$ as a direct limit of its coherent subsheaves.
For some index $i$, the induced map $(\shF_i)_{V,v'}\ra\O_V$ must be surjective.
Thus $\shF_i$ is a coherent sheaf on $\stG$ of weight $w=1$ locally containing invertible summands.
By Theorem \ref{gerbe equivalent}, the cohomology class $\alpha\in H^2(\stX,\GG_m)$ lies in the bigger Brauer group.
\qed

\medskip
The following example essentially due to Totaro (\cite{Totaro 2004}, Remark 1 in Introduction)
shows that the assumption on the diagonal  morphism $\Delta:\stX\ra\stX\times\stX$ in Theorem \ref{main result} is not
superfluous. Let $E$ be an elliptic curve over an algebraically closed ground field $k$,
and $\shL$ be an invertible sheaf on $E$ of degree zero, such that $\shL^{\otimes t}\neq \O_E$ for $t\neq 0$.
Consider the $\GG_{m,E}$-torsor $V=\Spec(\bigoplus_{t\in\ZZ}\shL^{\otimes t})$.
According to \cite{Serre 1975}, Chapter VII, \S 3.15, the torsor structure comes from
a unique extension of $k$-group schemes 
$0\ra\GG_m\ra V\ra E\ra 0$. From this we obtain a morphism of algebraic $k$-stacks
$BV\ra BE$, which sends a $V$-torsor to its associated $E$-torsor.
It follows that the morphism $BV\ra BE$ is a $\GG_{m,BE}$-gerbe.
Coherent sheaves on $BV$ correspond to   linear representations $V\ra\GL_k(n)$, $n\geq 0$.
Using that the scheme $\GL_k(n)$ is affine and $\Gamma(V,\O_V)=\bigoplus_{t\in\ZZ}\Gamma(E,\shL^{\otimes t})=k$,
we infer that every coherent sheaf on $BV$ is isomorphic to $\O_{BV}^{\oplus n}$.
In particular, there are no nonzero coherent sheaves of weight $w=1$.
Summing up, the algebraic $k$-stack $\stX=BE$ admits a $\GG_{m,\stX}$-gerbe
$\stG=BV$  whose cohomology class
does not lie in the bigger Brauer group.

\section{Appendix: Big-\'etale vs.\ lisse-\'etale cohomology}
\mylabel{appendix}

Let $\stX$ be an algebraic $S$-stack.
Then we have an inclusion   of sites $\Liset(\stX)\subset\Biget(\stX)$.
It obviously sends coverings to coverings,
whence the inclusion functor is continuous
by \cite{SGA 4a}, Expos\'e III, Proposition 1.6.
Hence for all big-\'etale sheaves $\shF$,
the induced presheaf $\shF_\liset =\shF|_{\Liset(\stX)}$ on the lisse-\'etale site   is  a sheaf.
Moreover, the induced restriction functor
$$
\stX_\biget\lra\stX_\liset,\quad \shF\longmapsto\shF_\liset
$$
commutes with all direct and inverse limits (by the formula for sheafification).
Consequently, the functors
$
\shF\mapsto H^i(\stX_\liset,\shF_\liset)
$
comprise a $\delta$-functor on the category of big-\'etale abelian sheaves.
By universality, the restriction map 
$\Gamma(\stX_\biget,\shF)\ra\Gamma(\stX_\liset,\shF_\liset)$
induces a natural transformation 
$$
H^i(\stX_\biget,\shF)\lra H^i(\stX_\liset,\shF_\liset)
$$
of $\delta$-functors. A priori, it is not clear that these canonical maps are
bijections, since there is no  map of topoi $u=(u_*,u^{-1})$ from the
big-\'etale to the lisse-\'etale topos,  with   $u_*(\shF)=\shF_\liset$.
The problem is as follows: By definition of   maps of topoi,
the functor $u^{-1}:\stX_\liset\ra\stX_\biget$ must be exact and left adjoint to $u_*$.
An adjoint indeed exists, and its values are the usual direct limits.
But as observed by Behrend and Gabber, the direct limits are not filtered,
whence the functor $u^{-1}$ is not left exact.
Compare the discussions in \cite{Behrend 2003}, Warning 4.42 and \cite{Olsson 2007}, Section 3.

The goal of this appendix is to establish that the canonical maps are nevertheless
bijections:

\begin{theorem}
\mylabel{cohomology bijective}
For all big-\'etale abelian sheaves $\shF$,
the canonical maps on cohomology groups $H^i(\stX_\biget,\shF)\ra H^i(\stX_\liset,\shF_\liset)$
are bijective.
\end{theorem}

To prove this statement, we shall   generalize it.
Recall that an abelian category $\shC$
satisfying Grothendieck's axiom $\text{AB5}$ (direct limits exists and are exact)
and containing a generator is  called a \emph{Grothendieck category}.
Typical examples are the category of modules over a ring,
or the category of abelian  sheaves  on a  site.
Recall that a \emph{generator}   is an object $U$ with
the following property: For every inclusion of objects $A\subset B$
there is a morphism $U\ra B $ not factoring over $A$.

\begin{lemma}
\mylabel{Grothendieck categories}
Let $F:\shC\ra\shD$ be an additive functor between Grothendieck
categories. Suppose that $F$ is   exact and commutes with all direct  sums,
and that any inclusion $B'\subset B$ of objects in $\shD$
is isomorphic to $F(A')\subset F(A)$ for some inclusion
$A'\subset A$ of objects in $\shC$.
Then for every object $A\in\shC$ there is an inclusion $A\subset I$
into an injective object $I\in\shC$ with the property that
$F(I)\in\shD$ is injective as well. 
\end{lemma}

\proof
We first observe that there is a generator $U\in\shC$ so that 
$F(U)\in\shD$ is a generator as well.
To see this, choose generators $U_1\in\shC$ and $V\in\shD$.
By assumption, there exists an object $U_2\in\shC$ with $F(U_2)=V$.
Then $U=U_1\oplus U_2$ and $F(U)=F(U_1)\oplus V$ are generators
in $\shC$ and $\shD$, respectively.

Next we recall  Grothendieck's   construction  of injective objects (\cite{Grothendieck 1957}, Theorem 1.10.1)
with a slight variant: 
Choose a family of injections $i_\alpha:U_\alpha\ra U$, $\alpha\in J$  
with the property that the set $i_\alpha(U_\alpha)\subset U$ runs
through the set of all subobjects. Here we allow repetitions, which do not occur in
the  original construction; this does not effect  the outcome of the construction, and
gives us a little extra freedom, which comes into play in the last paragraph.
Given an object $A\in \shC$,   let $J_A$ be the family
of all morphism $f_\beta:U_{\alpha(\beta)}\ra A$, $\beta\in J_A$
defined on some $U_{\alpha(\beta)}$, $\alpha(\beta)\in J$.
We now define another object $M(A)\in\shC$ and an injective morphism 
$A\ra M(A)$  
by the exact sequence
$$
\bigoplus_{\beta\in J_A} U_{\alpha(\beta)} \lra A \oplus 
\bigoplus_{\beta\in J_A} U\lra M(A)\lra 0.
$$
The map on the left is the canonical one:
$$
(u_\beta)\longmapsto (\sum f_\beta(u_\beta),(i_{\alpha(\beta)}(u_\beta))).
$$
Fix a cardinal $\mu$ that is at least as large as the cardinality
of the  set of all subobjects inside the generator $U$.
One defines, by transfinite induction, for any ordinal number
$\gamma\leq \mu$ a direct system 
$(M_\gamma(A))$  of objects
in $\shC$ as follows: 
$$
M_\gamma(A) =
\begin{cases}
A & \text{if $\gamma=0$,}\\
M(M_{\gamma'}(A)) & \text{if $\gamma=\gamma'+1$ is a successor ordinal,}\\
\dirlim_{\gamma'<\gamma} M_{\gamma'}(A) & \text{if $\gamma$ is a limit ordinal.}
\end{cases}
$$
Then the canonical maps $A=M_0(A)\ra M_\gamma(A)$ are injective,
and, by Grothendieck,  the object $M_\mu(A)$
is necessarily injective.
We recommend \cite{Krivine 1969} as a general reference for
ordinal and cardinal numbers.

The construction leads to the desired inclusion $A\subset I$
with $I$ and $F(I)$ injective:
By construction, $F(U)$ is a generator in $\shD$, the induced maps
$F(U_\alpha)\ra F(U)$ are injective, and their images run
through the set of all subobjects in $F(U)$, possibly with
repetitions.
The construction of $M_\gamma(A)$ uses only the formation
of direct limits, whence commutes with $F$, such that
$F(M_\gamma(A))\simeq M_\gamma(F(A))$.
So   both $I=M_\mu(A)$
and $F(I)=M_\mu(F(A))$ are   injective objects.
\qed

\begin{remark}
\mylabel{right inverse}
The condition that every inclusions $B'\subset B$ in $\shD$ is isomorphic
to $F(A')\subset F(A)$ for some $A'\subset A$
is satisfied  if there is an additive functor $G:\shD\ra\shC$
with the property that the composition $F\circ G$ is isomorphic to the identity $\id_\shD$:
Simply set $A=G(B)$ and let $A'\subset A$ be the image of
the induced morphism $G(B')\ra G(B)$.
\end{remark}
 
\emph{Proof of Theorem \ref{cohomology bijective}:}
We have to check that the canonical maps on cohomology
$H^i(\stX_\biget,\shF)\ra H^i(\stX_\liset,\shF_\liset)$ are
bijective in degree $i=0$, and that the
$\delta$-functor $\shF\mapsto H^i(\stX_\liset,\shF_\liset)$
is universal.

Indeed, that restriction map $\Gamma(\stX_\biget,\shF)\ra\Gamma(\stX_\liset,\shF_\liset)$ is bijective
for all set-valued big-\'etale sheaves $\shF$.
Recall that on sites like $\Biget(\stX)$ and $\Liset(\stX)$ that have no final object,
the global section functor is defined as a morphism set
$$
\Gamma(\stX_\biget,\shF)=\Mor(\underline{e},\shF),
$$
where $\underline{e}$ denotes the   sheaf whose values is constantly a 1-element set.
To proceed, choose a smooth surjection $u:U\ra\stX$ from some
algebraic space $U$, and let $\underline{U}$ be the big-\'etale sheaf
represented by $U$. According to \cite{SGA 4a}, Expos\'e II, Proposition 5.1,
the canonical map $\underline{U}\ra \underline{e}$
is an epimorphism in the topos $\stX_\biget$. By loc.\ cit.\, Expos\'e II,
Proposition 4.3, epimorphisms in topoi are effective and universal;
this simply means that the the sequence morphism sets
$\shF(\underline{e})\ra\shF(\underline{U})\rightrightarrows\shF(\underline{U}\times\underline{U})$
is exact.
This latter is nothing but
$$
\Gamma(\stX_\biget,\shF)\ra\Gamma(U,\shF)\rightrightarrows
\Gamma(U\times_\stX U,\shF).
$$
Using that $(U,u)$ is an object from the lisse-\'etale site, 
we easily   deduce  that the restriction map 
$\Gamma(\stX_\biget,\shF)\ra\Gamma(\stX_\liset,\shF_\liset)$ is bijective.

It remains to check that the  $\delta$-functor $\shF\mapsto H^i(\stX_\liset,\shF_\liset)$
is universal. According to \cite{Grothendieck 1957}, Proposition 2.2.1, it 
suffices to check that for every big-\'etale abelian sheaf $\shF$, there is
an inclusion $\shF\subset\shI$ into an injective sheaf so that the restriction 
$\shI_\liset$ is injective as well.

For  this  we construct a   functor
$(\Ab/\stX_\liset)\ra(\Ab/\stX_\biget)$, $\shG\mapsto\shG'$
as follows:
 Given an abelian lisse-\'etale   sheaf $\shG$, 
define an abelian big-\'etale    presheaf $\shG'$ by 
$$
\Gamma((U,u),\shG')=
\begin{cases}
\Gamma((U,u),\shG) & \text{if $u:U\ra \stX$ is smooth,}\\
0 & \text{otherwise.}
\end{cases}
$$
It is easy to see that $\shG'$ is a sheaf, and that
we have  a canonical map of sheaves $\shG\ra(\shG')_\liset$, which
is bijective. Note that the functor $\shG\mapsto\shG'$ resembles
the extension-by-zero functor for an open inclusion.
Remark \ref{right inverse} tells us that Lemma \ref{Grothendieck categories} applies to our situation. Whence
the desired inclusion $\shF\subset\shI$ with 
$\shI$ and $\shI_\liset$ injective exists.
\qed


\end{document}